\documentclass{article}
\usepackage{amsfonts}
\usepackage{amsmath}
\usepackage{harvard}

\newtheorem{theorem}{Theorem}

\newtheorem{corollary}[theorem]{Corollary}

\newtheorem{definition}[theorem]{Definition}

\newtheorem{lemma}[theorem]{Lemma}

\newtheorem{proposition}[theorem]{Proposition}
\newtheorem{remark}[theorem]{Remark}

\newenvironment{proof}[1][Proof]{\noindent\textbf{#1.} }{\ \rule{0.5em}{0.5em}}
\numberwithin{theorem}{section}
\numberwithin{equation}{section}
\input{tcilatex}

\begin{document}

\title{Jet Berwald-Riemann-Lagrange Geometrization for Affine Maps between
Finsler Manifolds}
\author{Mircea Neagu}
\date{}
\maketitle

\begin{abstract}
In this paper we introduce a natural definition for the affine maps between
two Finsler manifolds $(M,F)$ and $(N,\widetilde{F})$ and we give some
geometrical properties of these affine maps. Starting from the equations of
the affine maps, we construct a natural Berwald-Riemann-Lagrange geometry on
the 1-jet space $J^{1}(TM,N)$, in the sense of a Berwald nonlinear
connection $\Gamma _{\text{jet}}^{\text{b}}$, a Berwald $\Gamma _{\text{jet}%
}^{\text{b}}$-linear d-connection $B\Gamma _{\text{jet}}^{\text{b}}$,
together with its d-torsions and d-curvatures, which geometrically
characterizes the initial affine maps between Finsler manifolds.
\end{abstract}



\section{Introduction}

\hspace{5mm}It is well known that the harmonic maps between Riemannian
manifolds are defined as extremals of the energy functional. Because these
harmonic maps are very important in differential geometry and mathematical
physics, they were intensively studied by Eells and Lemaire [4].

By using the volume element induced on the projective sphere bundle $SM$ of
a Finsler manifold $(M,F)$, the harmonic maps between a Finsler manifold and
a Riemannian manifold were considered by Mo [7].

Recent studies on Finsler geometry led to the investigation of the
nondegenerate harmonic maps between two Finsler manifolds $(M,F)$ and $(N,%
\widetilde{F})$, as critical points of a natural energy functional on the
sphere bundle $SM.$ Thus, Shen and Zhang studied the variation formulas [11]
for harmonic maps between Finsler manifolds, and He and Shen established a
corresponding generalized Weitzenb\"{o}ck formula [5].

A general geometrical approach for harmonic maps between two generalized
Lagrange spaces $(M,g_{\alpha \beta }(t^{\gamma },s^{\gamma }))$ and $(N,%
\widetilde{g}_{ij}(x^{k},y^{k}))$ is given by the author of this paper in
[8].

In this work we investigate affine maps between Finsler manifolds, as
particular cases of nondegenerate harmonic maps. From a geometrical point of
view, we believe that our particular case of nondegenerate harmonic maps (we
refer to the affine maps) is not too restrictive one, because we consider
that there exist enough interesting geometrical results which characterize
the nondegenerate affine maps between Finsler manifolds. In this direction,
using the \textit{Riemann-Lagrange geometry on 1-jet spaces} recently
developed by the author of this paper in [9] and [10], we will show that the
equations of the nondegenerate affine maps between two Finsler manifolds $%
(M,F)$ and $(N,\widetilde{F})$ produce some natural d-torsions and
d-curvatures on the 1-jet space $J^{1}(TM,N)$, where $TM$ is the tangent
bundle of the smooth manifold $M.$

We would like to point out that the jet Riemann-Lagrange geometrical ideas
detailedly exposed in the works [9] and [10] were initially stated by Asanov
in the paper [2].

\section{Basic formulas on Finsler manifolds}

\hspace{5mm}Let us denote by $M$ a $p$-dimensional smooth manifold, which
induces on its tangent bundle $TM$ the local coordinates $(t^{\alpha
},s^{\alpha })$. Throughout this paper the greek indices $\alpha ,\beta
,\gamma ,...$ run from $1$ to $p$. Let us consider that the manifold $M$ is
endowed with a Finsler structure $F:TM\rightarrow \lbrack 0,\infty )$, such
that $(M,F)$ is a Finsler manifold.

The \textit{fundamental metrical d-tensor} of the Finsler manifold $(M,F)$
is defined on $TM\backslash \{0\}$ by%
\begin{equation*}
g_{\alpha \beta }(t^{\varepsilon },s^{\varepsilon })=\frac{1}{2}\frac{%
\partial ^{2}F^{2}}{\partial s^{\alpha }\partial s^{\beta }}.
\end{equation*}

\begin{remark}
Taking into account that $F^{2}$ is $2$-positive homogenous, we immediately
deduce, via the Euler theorem, that we have $F^{2}=g_{\alpha \beta
}s^{\alpha }s^{\beta }.$
\end{remark}

The fundamental metrical d-tensor $(g_{\alpha \beta })$ produces the \textit{%
Cartan d-tensor} of the Finsler manifold $(M,F)$, putting%
\begin{equation*}
C_{\alpha \beta \gamma }=\frac{1}{2}\frac{\partial g_{\alpha \beta }}{%
\partial s^{\gamma }}=\frac{1}{4}\frac{\partial ^{3}F^{2}}{\partial
s^{\alpha }\partial s^{\beta }\partial s^{\gamma }}.
\end{equation*}%
It is obvious that the Cartan d-tensor is totally symmetric in the indices $%
\alpha ,\beta $ and $\gamma .$ Moreover, because $g_{\alpha \beta }$'s are
positive homogenous of degree zero, the Euler theorem implies the equalities%
\begin{equation*}
C_{0\beta \gamma }=C_{\alpha 0\gamma }=C_{\alpha \beta 0}=0,
\end{equation*}%
where by the index $0$ we understand the contraction with $s^{\mu }$. For
instance, we have $C_{\alpha \beta 0}=C_{\alpha \beta \mu }s^{\mu }.$

Let us consider the \textit{formal Christoffel symbols} of the second kind%
\begin{equation}
\mathbf{\gamma }_{\alpha \beta }^{\mu }=\frac{g^{\mu \varepsilon }}{2}\left( 
\frac{\partial g_{\varepsilon \alpha }}{\partial t^{\beta }}+\frac{\partial
g_{\varepsilon \beta }}{\partial t^{\alpha }}-\frac{\partial g_{\alpha \beta
}}{\partial t^{\varepsilon }}\right) ,  \label{g}
\end{equation}%
where $(g^{\mu \varepsilon })$ denotes the inverse matrix of $(g_{\mu
\varepsilon })$. The formal Christoffel symbols produce on $TM\backslash
\{0\}$ the \textit{nonlinear Cartan connection} of the Finsler manifold $%
(M,F)$, taking (see [3], pp. 33)%
\begin{equation}
N_{\alpha }^{\beta }=\mathbf{\gamma }_{\alpha \varepsilon }^{\beta
}s^{\varepsilon }-C_{\alpha \varepsilon }^{\beta }\mathbf{\gamma }_{\mu \nu
}^{\varepsilon }s^{\mu }s^{\nu },  \label{CartanNLC}
\end{equation}%
where $C_{\alpha \varepsilon }^{\beta }=g^{\beta \lambda }C_{\lambda \alpha
\varepsilon }.$ An important geometrical concept in Finsler geometry is
given by the following notion:

\begin{definition}
A curve $c:[a,b]\rightarrow M$, locally expressed by $t^{\alpha }=t^{\alpha
}(t)$, where $t\in \lbrack a,b]$, is called an \textbf{autoparallel curve of
the nonlinear Cartan connection of the Finsler manifold} $(M,F)$ or,
briefly, an \textbf{autoparallel curve} on $(M,F)$, if and only if%
\begin{equation}
\frac{d^{2}t^{\alpha }}{dt^{2}}+N_{\beta }^{\alpha }\left( t^{\mu }(t),\frac{%
dt^{\mu }}{dt}\right) \frac{dt^{\beta }}{dt}=0,\text{ }\forall \text{ }%
\alpha =\overline{1,p}.  \label{auto}
\end{equation}
\end{definition}

Using the nonlinear Cartan connection $(N_{\beta }^{\alpha })$ we can
construct the \textit{generalized Christoffel symbols}%
\begin{equation*}
\Gamma _{\alpha \beta }^{\gamma }=\frac{g^{\gamma \mu }}{2}\left( \frac{%
\delta g_{\mu \alpha }}{\delta t^{\beta }}+\frac{\delta g_{\mu \beta }}{%
\delta t^{\alpha }}-\frac{\delta g_{\alpha \beta }}{\delta t^{\mu }}\right) ,
\end{equation*}%
where%
\begin{equation*}
\frac{\delta }{\delta t^{\alpha }}=\frac{\partial }{\partial t^{\alpha }}%
-N_{\alpha }^{\varepsilon }\frac{\partial }{\partial s^{\varepsilon }}.
\end{equation*}

\begin{remark}
The set of d-vector fields%
\begin{equation*}
\left\{ \frac{\delta }{\delta t^{\alpha }},\frac{\partial }{\partial
s^{\alpha }}\right\} \subset \mathcal{X}(TM\backslash \{0\})
\end{equation*}%
represents a basis in the set of vector fields on $TM\backslash \{0\}$,
which is called the \textbf{adapted basis} produced by the Finsler structure 
$F$. The transformation rules of the elements of the adapted basis are
tensorial ones.
\end{remark}

From a geometrical point of view, the generalized Christoffel symbols can be
regarded as, in Miron and Anastasiei's terminology [6, pp. 149], the \textit{%
adapted components} on $TM\backslash \{0\}$ of the \textit{distinguished
linear Rund connection}%
\begin{equation*}
R\Gamma =\left( N_{\alpha }^{\beta },\Gamma _{\alpha \beta }^{\gamma
},0\right) ,
\end{equation*}%
or, in Bao, Chern and Shen's terminology [3, pp. 39], as the coefficients of
the \textit{Chern connection} on the pulled-back tangent bundle $\pi ^{\ast
}TM\rightarrow TM\backslash \{0\}$, where $\pi :TM\backslash
\{0\}\rightarrow M$ is the canonical projection.

\begin{remark}
(i) For practitioners of Finsler geometry, Anastasiei pointed out in [1]
that the \textbf{Rund connection} and the \textbf{Chern connection}
coincide. In a such geometrical context, we underline that in this paper we
follow the terminology and the notations used by Miron and Anastasiei in [6].

(ii) It is important to note that, on the Finsler manifold $(M,F)$, the
formula%
\begin{equation*}
N_{\alpha }^{\beta }=\Gamma _{\alpha \gamma }^{\beta }s^{\gamma }
\end{equation*}%
is always true and very useful. Consequently, the \textbf{autoparallel curves%
} of the Finsler manifold $(M,F)$ are the solutions of the system of
differential equations%
\begin{equation}
\frac{d^{2}t^{\mu }}{dt^{2}}+\Gamma _{\alpha \beta }^{\mu }\left(
t^{\varepsilon }(t),\frac{dt^{\varepsilon }}{dt}\right) \frac{dt^{\alpha }}{%
dt}\frac{dt^{\beta }}{dt}=0,\text{ }\forall \text{ }\mu =\overline{1,p},
\label{autoG}
\end{equation}%
which is equivalent with the ODEs system of second order%
\begin{equation}
\frac{d^{2}t^{\mu }}{dt^{2}}+\mathbf{\gamma }_{\alpha \beta }^{\mu }\left(
t^{\varepsilon }(t),\frac{dt^{\varepsilon }}{dt}\right) \frac{dt^{\alpha }}{%
dt}\frac{dt^{\beta }}{dt}=0,\text{ }\forall \text{ }\mu =\overline{1,p}.
\label{autog}
\end{equation}

It follows that an \textbf{autoparallel curve} $c(t)=(t^{\alpha }(t))$
having a \textbf{natural parameter}, that is $F(c(t),dc/dt)=$ constant, is
equivalent with a \textbf{constant speed geodesic} on the Finsler manifold $%
(M,F).$ For more details, see [3, pp. 125-128] and [6, pp. 132-138].
\end{remark}

In the Finsler geometry literature, a famous linear d-connection on the
Finsler manifold $(M,F)$ is the \textit{Berwald connection}%
\begin{equation*}
B\Gamma =\left( N_{\alpha }^{\beta },B_{\alpha \beta }^{\gamma },0\right) ,
\end{equation*}%
whose \textit{adapted components} are defined by%
\begin{equation}
B_{\alpha \beta }^{\gamma }=\Gamma _{\alpha \beta }^{\gamma }+C_{\alpha
\beta |\mu }^{\gamma }s^{\mu }=\Gamma _{\alpha \beta }^{\gamma }+C_{\alpha
\beta |0}^{\gamma },  \label{B}
\end{equation}%
where%
\begin{equation*}
C_{\alpha \beta |\mu }^{\gamma }=\frac{\delta C_{\alpha \beta }^{\gamma }}{%
\delta t^{\mu }}+C_{\alpha \beta }^{\varepsilon }\Gamma _{\varepsilon \mu
}^{\gamma }-C_{\varepsilon \beta }^{\gamma }\Gamma _{\alpha \mu
}^{\varepsilon }-C_{\alpha \varepsilon }^{\gamma }\Gamma _{\beta \mu
}^{\varepsilon }
\end{equation*}%
represent the \textit{local horizontal covariant derivatives} produced by
the Rund connection $R\Gamma .$

\begin{remark}
(i) The local horizontal covariant derivatives produced by the Rund
connection $R\Gamma $ behave as differentiations, in the sense that they
work by the Leibniz rule upon the tensorial product of d-tensors and they
commute with the contractions of d-tensors.

(ii) The following special properties of the local horizontal covariant
derivatives associated to the Rund connection $R\Gamma $ hold good (see [6,
pp. 140]):%
\begin{equation*}
g_{\alpha \beta |\gamma }=0,\text{ }s_{\text{ \ }|\gamma }^{\alpha }=0,\text{
}F_{|\gamma }=0.
\end{equation*}

(iii) If we define the \textbf{spray} of the Finsler manifold $(M,F)$ by%
\begin{equation}
G^{\mu }=\frac{1}{2}\mathbf{\gamma }_{\alpha \beta }^{\mu }s^{\alpha
}s^{\beta }=\frac{1}{2}\Gamma _{\alpha \beta }^{\mu }s^{\alpha }s^{\beta },
\label{G}
\end{equation}%
then, by a direct calculation, we find the following relations:%
\begin{equation}
2G^{\gamma }=N_{\alpha }^{\gamma }s^{\alpha },\text{ }N_{\alpha }^{\gamma }=%
\frac{\partial G^{\gamma }}{\partial s^{\alpha }}\text{ and }B_{\alpha \beta
}^{\gamma }=\frac{\partial N_{\alpha }^{\gamma }}{\partial s^{\beta }}=\frac{%
\partial ^{2}G^{\gamma }}{\partial s^{\alpha }\partial s^{\beta }}.
\label{BG}
\end{equation}

(iv) The equations of the autoparallel curves of the Finsler manifold $(M,F)$
have the following simple form:%
\begin{equation*}
\frac{d^{2}t^{\alpha }}{dt^{2}}+2G^{\alpha }\left( t^{\mu }(t),\frac{dt^{\mu
}}{dt}\right) =0,\text{ }\forall \text{ }\alpha =\overline{1,p}.
\end{equation*}%
Moreover, following the geometrical ideas from [6, pp. 132-136], it is
important to note that, in fact, the autoparallel curves of the Finsler
manifold $(M,F)$ coincide exactly with the extremal curves of the integral
action of the absolute energy%
\begin{equation*}
\mathbf{E}(c)=\int_{a}^{b}F^{2}\left( t^{\alpha }(t),\frac{dt^{\alpha }}{dt}%
\right) dt.
\end{equation*}%
In this context, taking into account that $F^{2}=g_{\alpha \beta }s^{\alpha
}s^{\beta }$, a direct calculation of the Euler-Lagrange equations of the
preceding energy functional leads us to the formula%
\begin{equation}
G^{\gamma }(t^{\varepsilon },s^{\varepsilon })=\frac{g^{\gamma \mu }}{4}%
\left[ \frac{\partial ^{2}F^{2}}{\partial s^{\mu }\partial t^{\nu }}s^{\nu }-%
\frac{\partial F^{2}}{\partial t^{\mu }}\right] .  \label{sprayE}
\end{equation}
\end{remark}

From a Finsler geometrical point of view, we point out that the Berwald
connection $B\Gamma $ is characterized in the adapted basis $\{\delta
/\delta t^{\alpha },\partial /\partial s^{\alpha }\}$ by \textit{one} 
\textit{local torsion} d-tensor%
\begin{equation}
^{\mathbf{b}}\mathcal{R}_{\beta \gamma }^{\alpha }=\frac{\delta N_{\beta
}^{\alpha }}{\delta t^{\gamma }}-\frac{\delta N_{\gamma }^{\alpha }}{\delta
t^{\beta }}  \label{BT}
\end{equation}%
and \textit{two essential local curvature} d-tensors%
\begin{equation}
\begin{array}{l}
^{\mathbf{b}}\mathcal{R}_{\beta \gamma \varepsilon }^{\alpha }=\dfrac{\delta
B_{\beta \gamma }^{\alpha }}{\delta t^{\varepsilon }}-\dfrac{\delta B_{\beta
\varepsilon }^{\alpha }}{\delta t^{\gamma }}+B_{\beta \gamma }^{\mu }B_{\mu
\varepsilon }^{\alpha }-B_{\beta \varepsilon }^{\mu }B_{\mu \gamma }^{\alpha
}\medskip \\ 
^{\mathbf{b}}\mathcal{P}_{\beta \gamma \varepsilon }^{\text{ }\alpha }=%
\dfrac{\partial B_{\beta \gamma }^{\alpha }}{\partial s^{\varepsilon }}.%
\end{array}
\label{BC}
\end{equation}%
For more details, the reader is invited to compare the book [6, pp. 48, 122,
149] with the book [3, pp. 52, 67].

\begin{remark}
Taking into account the relations (\ref{BG}), note that the Berwald
curvature d-tensor $^{\mathbf{b}}\mathcal{P}_{\beta \gamma \varepsilon }^{%
\text{ }\alpha }$ is totally symmetric in the indices $\beta ,\gamma $ and $%
\varepsilon $.
\end{remark}

\section{Affine maps between Finsler manifolds}

\hspace{5mm}Let $(M,F)$ and $(N,\widetilde{F})$ be two Finsler manifolds,
where the dimension of $N$ is $n$, and let $\varphi :(M,F)\rightarrow (N,%
\widetilde{F})$ be a smooth map which is \textit{nondegenerate}, that is its
induced tangent map verifies the condition $Ker(d\varphi )=\{0\}.$

\begin{remark}
(i) We suppose that the tangent bundle $TN$ has the local coordinates $%
(x^{i},y^{i})$, where $i=\overline{1,n}.$ Moreover, throughout this paper,
we suppose that the latin indices $i,j,k,...$ run from $1$ to $n$.

(ii) On the source Finsler manifold $(M,F)$ we will use the notations and
indices from Section 2, and on the target Finsler manifold $(N,\widetilde{F}%
) $ we will denote the same geometrical entities by the same letters, but
with tilde and corresponding indices.
\end{remark}

Let us suppose that the nondegenerate smooth map $\varphi $ is locally
expressed by $\varphi ^{i}=\varphi ^{i}(t^{\mu })$ and let us introduce the
notations%
\begin{equation*}
\varphi _{\alpha }^{i}=\frac{\partial \varphi ^{i}}{\partial t^{\alpha }}%
\text{ and }\varphi _{\alpha \beta }^{i}=\frac{\partial ^{2}\varphi ^{i}}{%
\partial t^{\alpha }\partial t^{\beta }}.
\end{equation*}%
In this geometrical context, we introduce the following concept:

\begin{definition}
The nondegenerate smooth map $\varphi :(M,F)\rightarrow (N,\widetilde{F})$
is called an \textbf{affine map between the Finsler manifolds }$(M,F)$%
\textbf{\ and }$(N,\widetilde{F})$ if and only if%
\begin{equation}
\varphi _{\alpha \beta }^{i}-B_{\alpha \beta }^{\gamma }\varphi _{\gamma
}^{i}+\widetilde{B}_{jk}^{i}\varphi _{\alpha }^{j}\varphi _{\beta }^{k}=0,%
\text{ }\forall \text{ }\alpha ,\beta =\overline{1,p},\text{ }\forall \text{ 
}i=\overline{1,n},  \label{Af}
\end{equation}%
where%
\begin{equation*}
B_{\alpha \beta }^{\gamma }=B_{\alpha \beta }^{\gamma }(t^{\mu },s^{\mu })%
\text{ and }\widetilde{B}_{jk}^{i}=\widetilde{B}_{jk}^{i}(\varphi
^{l}(t^{\mu }),\varphi _{\varepsilon }^{l}(t^{\mu })s^{\varepsilon })
\end{equation*}%
are the adapted components of the Berwald connection on $(M,F)$ and $(N,%
\widetilde{F})$, respectively.
\end{definition}

\begin{remark}
(i) If the target Finsler manifold $(N,\widetilde{F})$ is a Riemannian one,
then we have $\widetilde{C}_{jk}^{i}=0$ and $\widetilde{B}_{jk}^{i}=%
\widetilde{\gamma }_{jk}^{i}(\varphi ^{l}(t^{\mu }))$. In this case, the
assumption on the nondegeneration of $\varphi $ is not necessary in our
definition.

(ii) If $(M,F)$ and $(N,\widetilde{F})$ are both Riemannian manifolds, then
we recover the classical definition of the affine maps between two
Riemannian manifolds.
\end{remark}

In order to give some geometrical examples of affine maps between two
Finsler manifolds, let us introduce the following definition:

\begin{definition}
A smooth map $\varphi :(M,F)\rightarrow (N,\widetilde{F})$ is called an 
\textbf{isometry between the Finsler manifolds }$(M,F)$\textbf{\ and }$(N,%
\widetilde{F})$ or, briefly, \textbf{Finsler isometry}, if it satisfies the
conditions:\medskip

(i) $\varphi $ is a difeomorphism;\medskip

(ii) $F(t,s)=\widetilde{F}(\varphi (t),d\varphi (s)),$ $\forall $ $(t,s)\in
TM\backslash \{0\}$.
\end{definition}

\begin{remark}
If the smooth map $\varphi :M\rightarrow N$ is a difeomorphism, then its
induced tangent map $d\varphi _{t}:T_{t}M\rightarrow T_{\varphi (t)}N$, $%
\forall $ $t\in M$, is an isomorphism of vector spaces. It follows that we
have $p=n$ and $\det (\varphi _{\alpha }^{i})\neq 0$, that is $\varphi $ is
a nondegenerate map.
\end{remark}

\begin{theorem}
Any Finsler isometry $\varphi :(M,F)\rightarrow (N,\widetilde{F})$ is an
affine map between the Finsler manifolds $(M,F)$ and $(N,\widetilde{F})$.
\end{theorem}

\begin{proof}
Let $\varphi :(M,F)\rightarrow (N,\widetilde{F})$ be a Finsler isometry and
let $(\psi _{i}^{\alpha })$ be the inverse of the matrix $(\varphi _{\alpha
}^{i})$. It follows that we have $\psi _{i}^{\alpha }\varphi _{\alpha
}^{j}=\delta _{i}^{j}$ and $\psi _{i}^{\alpha }\varphi _{\beta }^{i}=\delta
_{\beta }^{\alpha }.$

Because $\varphi $ is a Finsler isometry, a direct calculation leads us to
the relations%
\begin{equation}
g_{\alpha \beta }=\widetilde{g}_{ij}\varphi _{\alpha }^{i}\varphi _{\beta
}^{j}\text{ and }g^{\alpha \beta }=\widetilde{g}^{ij}\psi _{i}^{\alpha }\psi
_{j}^{\beta },\text{ }\forall \text{ }\alpha ,\beta =\overline{1,p}.
\label{M}
\end{equation}

Using formula (\ref{g}) and the relations (\ref{M}), a new direct
calculation gives us the relations%
\begin{eqnarray}
\mathbf{\gamma }_{\alpha \beta }^{\mu } &=&\widetilde{\mathbf{\gamma }}%
_{ij}^{m}\varphi _{\alpha }^{i}\varphi _{\beta }^{j}\psi _{m}^{\mu }+%
\widetilde{C}_{ij}^{m}\left( \varphi _{\alpha }^{i}\varphi _{\beta \nu
}^{j}+\varphi _{\beta }^{i}\varphi _{\alpha \nu }^{j}\right) \psi _{m}^{\mu
}s^{\nu }-  \label{gg} \\
&&-\widetilde{C}_{ijm}g^{\mu \varepsilon }\varphi _{\alpha }^{i}\varphi
_{\beta }^{j}\varphi _{\varepsilon \nu }^{m}s^{\nu }+\varphi _{\alpha \beta
}^{m}\psi _{m}^{\mu },\text{ }\forall \text{ }\alpha ,\beta ,\mu =\overline{%
1,p}.  \notag
\end{eqnarray}%
Contracting the relations (\ref{gg}) with $s^{\alpha }$ and $s^{\beta },$
the formula (\ref{G}) implies the equalities%
\begin{equation}
2G^{\gamma }=2\widetilde{G}^{m}\psi _{m}^{\gamma }+\varphi _{\alpha \beta
}^{m}\psi _{m}^{\gamma }s^{\alpha }s^{\beta },\text{ }\forall \text{ }\gamma
=\overline{1,p},  \label{GG}
\end{equation}%
where, taking into account that we have $\widetilde{C}_{ijk}=\widetilde{C}%
_{ijk}(\varphi ^{l},\varphi _{\varepsilon }^{l}s^{\varepsilon }),$ we used
the relation 
\begin{equation*}
\widetilde{C}_{ijk}\varphi _{\varepsilon }^{k}s^{\varepsilon }=0.
\end{equation*}%
In the sequel, a double differentiation in (\ref{GG}), together with the
relation (\ref{BG}), imply the equalities%
\begin{equation}
B_{\alpha \beta }^{\gamma }=\widetilde{B}_{jk}^{m}\varphi _{\alpha
}^{j}\varphi _{\beta }^{k}\psi _{m}^{\gamma }+\varphi _{\alpha \beta
}^{m}\psi _{m}^{\gamma },\text{ }\forall \text{ }\alpha ,\beta ,\gamma =%
\overline{1,p}.  \label{BB}
\end{equation}

It is obvious now that the equalities (\ref{BB}) imply the equalities (\ref%
{Af}), which represent the equations of the affine maps. In conclusion, $%
\varphi $ is an affine map between the Finsler manifolds $(M,F)$ and $(N,%
\widetilde{F})$.
\end{proof}

\begin{corollary}
The identity map $\mathcal{I}:(N,\widetilde{F})\rightarrow (N,\widetilde{F})$
is an affine map.
\end{corollary}

Now, let us study the affinity of the identity map $\mathcal{I}$ when it
works with two different Finsler structures $F$ and $\widetilde{F}$ on the
manifold $N$.

\begin{proposition}
\label{Id}The identity map $\mathcal{I}:(N,F)\rightarrow (N,\widetilde{F})$
is an affine map if and only if 
\begin{equation}
G^{i}=\widetilde{G}^{i},\text{ }\forall \text{ }i=\overline{1,n},  \label{eg}
\end{equation}%
where $G^{i}$ and $\widetilde{G}^{i}$ are the spray coefficients of the
Finsler manifolds $(N,F)$ and $(N,\widetilde{F})$, respectively.
\end{proposition}

\begin{proof}
Let us suppose that we locally have%
\begin{equation*}
\mathcal{I}^{i}=\mathcal{I}^{i}(x^{k})=x^{i},\text{ }\forall \text{ }i=%
\overline{1,n}.
\end{equation*}%
Then, it immediately follows that we have%
\begin{equation*}
\mathcal{I}_{j}^{i}=\delta _{j}^{i}\text{ and }\mathcal{I}_{jk}^{i}=0.
\end{equation*}%
Consequently, $\mathcal{I}$ is an affine map between the Finsler manifolds $%
(N,F)$ and $(N,\widetilde{F})$ if and only if we have%
\begin{equation}
B_{jk}^{i}(x^{l},y^{l})=\widetilde{B}_{jk}^{i}(x^{l},y^{l})\Leftrightarrow 
\frac{\partial ^{2}G^{i}}{\partial y^{j}\partial y^{k}}=\frac{\partial ^{2}%
\widetilde{G}^{i}}{\partial y^{j}\partial y^{k}},\text{ }\forall \text{ }%
i,j,k=\overline{1,n},  \label{DG}
\end{equation}%
where%
\begin{equation*}
G^{i}=\frac{1}{2}\mathbf{\gamma }_{pq}^{i}y^{p}y^{q}\text{ and }\widetilde{G}%
^{i}=\frac{1}{2}\widetilde{\mathbf{\gamma }}_{pq}^{i}y^{p}y^{q}.
\end{equation*}

Taking into account that the spray coefficients are $2$-positive homogenous,
by contractions with $y^{j}$ and $y^{k}$, the equalities (\ref{DG}) imply
the equalities (\ref{eg}). Conversely, it is obvious that the equalities (%
\ref{eg}) imply the equalities (\ref{DG}). In conclusion, we obtain what we
were looking for.
\end{proof}

\begin{corollary}
The identity map $\mathcal{I}:(N,F)\rightarrow (N,\widetilde{F})$ is an
affine map if and only if%
\begin{equation*}
g^{ip}\left[ \frac{\partial ^{2}F^{2}}{\partial y^{p}\partial x^{q}}y^{q}-%
\frac{\partial F^{2}}{\partial x^{p}}\right] =\widetilde{g}^{ip}\left[ \frac{%
\partial ^{2}\widetilde{F}^{2}}{\partial y^{p}\partial x^{q}}y^{q}-\frac{%
\partial \widetilde{F}^{2}}{\partial x^{p}}\right] ,\text{ }\forall \text{ }%
i=\overline{1,n}.
\end{equation*}
\end{corollary}

\begin{proof}
The Corollary is an immediate consequence of the Proposition \ref{Id} and
the formula (\ref{sprayE}).
\end{proof}

\begin{corollary}
If $(N,h)$ is a flat Riemannian manifold and $(N,\widetilde{F})$ is a
locally Minkowski manifold, then the identity maps 
\begin{equation*}
\mathcal{I}:(N,h)\rightarrow (N,\widetilde{F})
\end{equation*}%
and 
\begin{equation*}
\mathcal{I}:(N,\widetilde{F})\rightarrow (N,h)
\end{equation*}%
are affine maps.
\end{corollary}

\begin{proof}
If $\widetilde{F}$ is a local Minkowski structure, then there exists a
system of local coordinates $(x^{i})$ such that $\widetilde{F}(x,y)=%
\widetilde{F}(y)$, that is $\widetilde{g}_{ij}(x,y)=\widetilde{g}_{ij}(y)$.
It follows that we have $\widetilde{\mathbf{\gamma }}_{ij}^{k}=0$, that is $%
\widetilde{G}^{k}=0$, $\forall $ $k=\overline{1,n}.$

On the other hand, if the Riemannian metric $h$ is flat, then we also have $%
\mathbf{\gamma }_{ij}^{k}=0$, that is $G^{k}=0$, $\forall $ $k=\overline{1,n}%
.$
\end{proof}

\section{Some geometrical properties of the affine maps between Finsler
manifolds}

\hspace{5mm}In our geometrical context, let us consider the particular case
when our source Finsler manifold is the Euclidian manifold%
\begin{equation*}
(M,F)=(\mathbb{R},F(t,s)=|s|).
\end{equation*}

Then, we can prove the following result:

\begin{proposition}
Any affine map $c:(\mathbb{R},F)\rightarrow (N,\widetilde{F})$ is an
autoparallel curve of the Finsler manifold $(N,\widetilde{F}).$
\end{proposition}

\begin{proof}
If $M=\mathbb{R}$ is regarded as the Euclidian manifold $(\mathbb{R},1)$,
then the equations (\ref{Af}) of the affine maps become%
\begin{equation}
\frac{d^{2}c^{i}}{dt^{2}}+\widetilde{B}_{jk}^{i}\left( c^{l}(t),\frac{dc^{l}%
}{dt}s\right) \frac{dc^{j}}{dt}\frac{dc^{k}}{dt}=0,\text{ }\forall \text{ }i=%
\overline{1,n},\text{ }\forall \text{ }s\in \mathbb{R}^{\ast },  \label{AfRs}
\end{equation}%
where $c(t)=(c^{i}(t))$ is an affine map. Taking into account that $%
\widetilde{B}_{jk}^{i}$'s are $0$-positive homogenous, we deduce from
equations (\ref{AfRs}) that an affine map $c(t)$ must verify the equations%
\begin{equation}
\frac{d^{2}c^{i}}{dt^{2}}+\widetilde{B}_{jk}^{i}\left( c^{l}(t),\frac{dc^{l}%
}{dt}\right) \frac{dc^{j}}{dt}\frac{dc^{k}}{dt}=0,\text{ }\forall \text{ }i=%
\overline{1,n}.  \label{AfR}
\end{equation}%
Now, using the formula (\ref{B}) and the fact that 
\begin{equation*}
\widetilde{C}_{jk|0}^{i}(dc^{j}/dt)=0,
\end{equation*}%
it follows that the equations (\ref{AfR}) become exactly the equations (\ref%
{autoG}) of the autoparallel curves of the Finsler manifold $(N,\widetilde{F}%
).$
\end{proof}

In order to obtain a geometrical result which characterizes the affine maps
between Finsler manifolds, let us prove the following helpful statement:

\begin{lemma}
\label{Lema} Let $u_{\alpha \beta }:TM\backslash \{0\}\rightarrow \mathbb{R}$%
, where $\alpha ,\beta =\overline{1,p},$ be a family of smooth maps which
have the following four properties:\medskip

(i) $u_{\alpha \beta }=u_{\beta \alpha }$;\medskip

(ii) $u_{\alpha \beta }$'s are $0$-positive homogenous;\medskip

(iii) $\dfrac{\partial u_{\alpha \beta }}{\partial s^{\gamma }}$ is totally
symmetric in $\alpha ,\beta ,\gamma $;\medskip

(iv) $u_{\alpha \beta }s^{\alpha }s^{\beta }=0$.\medskip

Then, we have 
\begin{equation*}
u_{\alpha \beta }=0,\text{ }\forall \text{ }\alpha ,\beta =\overline{1,p}.
\end{equation*}
\end{lemma}

\begin{proof}
Differentiating (iv) with respect to $s^{\gamma }$ and using (i), we obtain
the equalities%
\begin{equation}
\frac{\partial u_{\alpha \beta }}{\partial s^{\gamma }}s^{\alpha }s^{\beta
}+2u_{\alpha \gamma }s^{\alpha }=0,\text{ }\forall \text{ }\gamma =\overline{%
1,p}.  \label{E1}
\end{equation}%
Using now (iii) and (ii), we immediately deduce that%
\begin{equation*}
\frac{\partial u_{\alpha \beta }}{\partial s^{\gamma }}s^{\alpha }s^{\beta }=%
\frac{\partial u_{\beta \gamma }}{\partial s^{\alpha }}s^{\alpha }s^{\beta
}=0,\text{ }\forall \text{ }\gamma =\overline{1,p}.
\end{equation*}%
Consequently, the equalities (\ref{E1}) become 
\begin{equation}
u_{\alpha \gamma }s^{\alpha }=0,\text{ }\forall \text{ }\gamma =\overline{1,p%
}.  \label{E2}
\end{equation}

Applying the same procedure to the equalities (\ref{E2}), we find that%
\begin{equation*}
u_{\gamma \varepsilon }=0,\text{ }\forall \text{ }\gamma ,\varepsilon =%
\overline{1,p}.
\end{equation*}
\end{proof}

\begin{theorem}[characterization of affine maps]
A nondegenerate smooth map $\varphi :(M,F)\rightarrow (N,\widetilde{F})$ is
an affine map between the Finsler manifolds $(M,F)$ and $(N,\widetilde{F})$
if and only if the map $\varphi $ carries autoparallel curves from $(M,F)$
into autoparallel curves on $(N,\widetilde{F})$.
\end{theorem}

\begin{proof}
Let $c(t)=(t^{\alpha }(t))$ be an autoparallel curve on the Finsler manifold 
$(M,F)$, that is it verifies the equations (\ref{autoG}).

Let us consider that the curve $\widetilde{c}(t)=(\varphi \circ c)(t)$ is
locally expressed by the components $x^{i}(t)=\varphi ^{i}(t^{\alpha }(t))$.
Then, differentiating by $t$, we immediately find for each $i\in \{1,...,n\}$
the equalities%
\begin{equation}
\frac{dx^{i}}{dt}=\varphi _{\alpha }^{i}\frac{dt^{\alpha }}{dt}\text{ and }%
\frac{d^{2}x^{i}}{dt^{2}}=\varphi _{\alpha \beta }^{i}\frac{dt^{\alpha }}{dt}%
\frac{dt^{\beta }}{dt}+\varphi _{\mu }^{i}\frac{d^{2}t^{\mu }}{dt^{2}}.
\label{d}
\end{equation}%
Using the equalities (\ref{autoG}), it follows that we have%
\begin{eqnarray}
\frac{d^{2}x^{i}}{dt^{2}}+\widetilde{\Gamma }_{jk}^{i}\frac{dx^{j}}{dt}\frac{%
dx^{k}}{dt} &=&\varphi _{\mu }^{i}\frac{d^{2}t^{\mu }}{dt^{2}}+\left(
\varphi _{\alpha \beta }^{i}+\widetilde{\Gamma }_{jk}^{i}\varphi _{\alpha
}^{j}\varphi _{\beta }^{k}\right) \frac{dt^{\alpha }}{dt}\frac{dt^{\beta }}{%
dt}=  \label{e} \\
&=&\left( \varphi _{\alpha \beta }^{i}-\Gamma _{\alpha \beta }^{\gamma
}\varphi _{\gamma }^{i}+\widetilde{\Gamma }_{jk}^{i}\varphi _{\alpha
}^{j}\varphi _{\beta }^{k}\right) \frac{dt^{\alpha }}{dt}\frac{dt^{\beta }}{%
dt}.  \notag
\end{eqnarray}

$"\Longrightarrow "$ If $\varphi $ is an affine map, then it verifies the
equations%
\begin{equation}
\varphi _{\alpha \beta }^{i}-\left( \Gamma _{\alpha \beta }^{\gamma
}+C_{\alpha \beta |0}^{\gamma }\right) \varphi _{\gamma }^{i}+\left( 
\widetilde{\Gamma }_{jk}^{i}+\widetilde{C}_{jk|0}^{i}\right) \varphi
_{\alpha }^{j}\varphi _{\beta }^{k}=0,  \label{Afdezv}
\end{equation}%
for any $\alpha ,\beta =\overline{1,p}$ and $i=\overline{1,n}.$ Contracting
the equations (\ref{Afdezv}) with $dt^{\alpha }/dt$ and $dt^{\beta }/dt$ and
taking into account the relations%
\begin{equation}
C_{\alpha \beta |0}^{\gamma }\frac{dt^{\alpha }}{dt}\frac{dt^{\beta }}{dt}=0,%
\text{ }\forall \text{ }\gamma =\overline{1,p},  \label{f1}
\end{equation}%
and%
\begin{equation}
\widetilde{C}_{jk|0}^{i}\varphi _{\alpha }^{j}\varphi _{\beta }^{k}\frac{%
dt^{\alpha }}{dt}\frac{dt^{\beta }}{dt}=0,\text{ }\forall \text{ }i=%
\overline{1,n},  \label{f2}
\end{equation}%
we deduce from equalities (\ref{e}) that we have%
\begin{equation*}
\frac{d^{2}x^{i}}{dt^{2}}+\widetilde{\Gamma }_{jk}^{i}\frac{dx^{j}}{dt}\frac{%
dx^{k}}{dt}=0,\text{ }\forall \text{ }i=\overline{1,n}.
\end{equation*}%
This is exactly what we were looking for.

$"\Longleftarrow "$ Conversely, if $\varphi $ carries autoparallel curves
from $(M,F)$ into autoparallel curves on $(N,\widetilde{F})$, then the
equalities (\ref{e}) imply the relations%
\begin{equation*}
\left( \varphi _{\alpha \beta }^{i}-\Gamma _{\alpha \beta }^{\gamma }\varphi
_{\gamma }^{i}+\widetilde{\Gamma }_{jk}^{i}\varphi _{\alpha }^{j}\varphi
_{\beta }^{k}\right) \frac{dt^{\alpha }}{dt}\frac{dt^{\beta }}{dt}=0,\text{ }%
\forall \text{ }i=\overline{1,n},
\end{equation*}%
for an arbitrary d-tensor field $s^{\alpha }=dt^{\alpha }/dt$. Obviously,
the relations (\ref{f1}) and (\ref{f2}) lead us to the equalities%
\begin{equation}
\left[ \varphi _{\alpha \beta }^{i}-\left( \Gamma _{\alpha \beta }^{\gamma
}+C_{\alpha \beta |0}^{\gamma }\right) \varphi _{\gamma }^{i}+\left( 
\widetilde{\Gamma }_{jk}^{i}+\widetilde{C}_{jk|0}^{i}\right) \varphi
_{\alpha }^{j}\varphi _{\beta }^{k}\right] s^{\alpha }s^{\beta }=0,
\label{Afu}
\end{equation}%
for any $i=\overline{1,n}.$ Denoting now the square parentheses from the
left side of the equalities (\ref{Afu}) with $u_{\alpha \beta }^{i},$ via
the relations (\ref{B}) and (\ref{BG}), we remark that we can apply the
Lemma \ref{Lema} to $u_{\alpha \beta }^{i},$ for any $i=\overline{1,n}$. In
conclusion, we obtain 
\begin{equation*}
u_{\alpha \beta }^{i}=0,\text{ }\forall \text{ }\alpha ,\beta =\overline{1,p}%
,\text{ }\forall \text{ }i=\overline{1,n},
\end{equation*}%
that is the nondegenerate map $\varphi $ is an affine map between the
Finsler manifolds $(M,F)$ and $(N,\widetilde{F})$.
\end{proof}

In the sequel, we will show that the affine maps are only particular cases
of \textit{harmonic maps between Finsler manifolds}, studied by Mo [7] or
Shen and Zhang [11].

\begin{proposition}
\label{AfHarm} Supposing that $M$ is a compact oriented smooth manifold
without boundary, then any affine map $\varphi :(M,F)\rightarrow (N,%
\widetilde{F})$ is a harmonic map between the Finsler manifolds $(M,F)$ and $%
(N,\widetilde{F})$, with vanishing tension field.
\end{proposition}

\begin{proof}
Following the geometrical ideas developed in [11], a particular case of
harmonic maps between the Finsler manifolds $(M,F)$ and $(N,\widetilde{F})$
is when the tension field of the nondegenerate smooth map $\varphi $
vanishes identically. But, the \textit{tension field} of $\varphi $ is given
by the components (see [11, pp. 45, (2.25)])%
\begin{eqnarray}
\tau ^{i}(\varphi ) &=&g^{\alpha \beta }\left\{ \varphi _{\alpha \beta
}^{i}-B_{\alpha \beta }^{\gamma }\varphi _{\gamma }^{i}+\widetilde{B}%
_{jk}^{i}\varphi _{\alpha }^{j}\varphi _{\beta }^{k}\right\} +  \label{TF} \\
&&+4g^{\alpha \beta }\widetilde{C}_{jk}^{i}\varphi _{\alpha }^{k}\left\{
\varphi _{\beta \gamma }^{j}-\Gamma _{\beta \gamma }^{\mu }\varphi _{\mu
}^{j}+\widetilde{\Gamma }_{pq}^{j}\varphi _{\beta }^{p}\varphi _{\gamma
}^{q}\right\} s^{\gamma }+  \notag \\
&&+g^{\alpha \beta }\widetilde{C}_{jkl}^{i}\varphi _{\alpha }^{k}\varphi
_{\beta }^{l}\left\{ \varphi _{\gamma \varepsilon }^{j}-\Gamma _{\gamma
\varepsilon }^{\mu }\varphi _{\mu }^{j}+\widetilde{\Gamma }_{pq}^{j}\varphi
_{\gamma }^{p}\varphi _{\varepsilon }^{q}\right\} s^{\gamma }s^{\varepsilon
},  \notag
\end{eqnarray}%
where%
\begin{equation*}
\widetilde{C}_{jkl}^{i}(x,y)=\widetilde{g}^{im}\frac{\partial \widetilde{C}%
_{jkl}}{\partial y^{m}}
\end{equation*}%
and, we underline that, in the expression (\ref{TF}), the Finsler
geometrical entities on $(N,\widetilde{F})$ are computed in 
\begin{equation*}
(x^{l},y^{l})=(\varphi ^{l}(t^{\mu }),\varphi _{\mu }^{l}s^{\mu })\in
TN\backslash \{0\}.
\end{equation*}%
Because we have%
\begin{equation*}
C_{\beta \gamma |0}^{\alpha }s^{\gamma }=\widetilde{C}_{jk|0}^{i}y^{k}=0,
\end{equation*}%
the formula (\ref{B}) implies that the components (\ref{TF}) of the tension
field of the nondegenerate smooth map $\varphi $ take the simpler form%
\begin{equation}
\tau ^{i}(\varphi )=g^{\alpha \beta }\tau _{\alpha \beta }^{i}+4g^{\alpha
\beta }\widetilde{C}_{jk}^{i}\varphi _{\alpha }^{k}\tau _{\beta \gamma
}^{j}s^{\gamma }+g^{\alpha \beta }\widetilde{C}_{jkl}^{i}\varphi _{\alpha
}^{k}\varphi _{\beta }^{l}\tau _{\gamma \varepsilon }^{j}s^{\gamma
}s^{\varepsilon },  \label{TFtau}
\end{equation}%
where%
\begin{equation*}
\tau _{\alpha \beta }^{i}=\varphi _{\alpha \beta }^{i}-B_{\alpha \beta
}^{\gamma }\varphi _{\gamma }^{i}+\widetilde{B}_{jk}^{i}\varphi _{\alpha
}^{j}\varphi _{\beta }^{k},\text{ }\forall \text{ }\alpha ,\beta =\overline{%
1,p},\text{ }\forall \text{ }i=\overline{1,n}.
\end{equation*}

It is obvious now that if the nondegenerate smooth map $\varphi $ is an
affine map, then its tension field vanishes identically. In conclusion, we
obtained what we were looking for.
\end{proof}

\begin{remark}
In the Proposition \ref{AfHarm} the assumptions upon the source manifold $M$
were imposed by the good definition of the harmonic maps between two Finsler
manifolds (see [11]). We point out yet that our definition of the affine
maps between Finsler manifolds did not require these assumptions.
\end{remark}

\begin{corollary}
Under the assumptions of Proposition \ref{AfHarm}, any nondegenerate smooth
map $\varphi :(M,F)\rightarrow (N,\widetilde{F})$, which carries
autoparallel curves from $(M,F)$ into autoparallel curves on $(N,\widetilde{F%
})$, is a harmonic map between the Finsler manifolds $(M,F)$ and $(N,%
\widetilde{F})$, with vanishing tension field.
\end{corollary}

\section{A Berwald-Riemann-Lagrange geometrization on the 1-jet space $%
J^{1}(TM,N)$, produced by the equations of the affine maps between two
Fins\-ler manifolds}

\hspace{5mm}The aim of this Section is to associate to the affine maps
equations (\ref{Af}) some geometrical objects which may characterize the
affine maps between two Finsler manifolds. Taking into account that starting
from the Asanov's geometrical ideas [2] the author of this paper has
recently elaborated a \textit{Riemann-Lagrange geometry on 1-jet spaces}
[9], in the sense of d-connections, d-torsions and d-curvatures, our
geometrical construction is made on the 1-jet space $J^{1}(TM,N)$. For a
more clear exposition of our jet geometrical ideas, it is important as the
tangent bundle $TM$ to be regarded as having the local coordinates $%
(t^{\alpha },s^{a})$, where the indices $\alpha ,\beta ,\gamma ,...$ and $%
a,b,c,...$ have the same range: $1,2,...,p.$

\begin{remark}
(i) In our geometrical context, it is more convenient to have two kinds of
indices, in order to mark the distinct elements of the \textbf{adapted basis}%
\begin{equation}
\left\{ \frac{\delta }{\delta t^{\alpha }}=\frac{\partial }{\partial
t^{\alpha }}-N_{\alpha }^{b}\frac{\partial }{\partial s^{b}},\frac{\partial 
}{\partial s^{a}}\right\} \subset \mathcal{X}(TM\backslash \{0\})  \label{AB}
\end{equation}%
provided by the canonical nonlinear Cartan connection (\ref{CartanNLC}) of
the Finsler structure $F:TM\rightarrow \mathbb{R}_{+}$ on the smooth
manifold $M$.

(ii) Note that we will use the formal notation $A=(\alpha ,a)$ for an index
which run two times from $1$ to $p$, namely, firstly by $\alpha $,
corresponding to the coordinates $t^{\alpha }$ or, equivalent, to the 
\textbf{horizontal d-vector fields} $\delta /\delta t^{\alpha }$, and after
that by $a$, corresponding to the coordinates $s^{a}$ or, equivalent, to the 
\textbf{vertical d-vector fields} $\partial /\partial s^{a}.$ Moreover,
throughout this paper, the capital latin letters $A,B,C,...$ will denote
indices like the previous one.
\end{remark}

Using the notation $\left( T^{A}\right) =\left( t^{\alpha },s^{a}\right) $
for the local coordinates on the tangent bundle $TM$, we recall that the
coordinates on the 1-jet space $J^{1}(TM,N)$ are $\left(
T^{A},x^{i},X_{A}^{i}\right) $, where the coordinates $\left(
X_{A}^{i}\right) =\left( x_{\alpha }^{i},y_{a}^{i}\right) $ have the meaning
of the \textit{partial derivatives} of the functions $x^{i}$ with respect to 
$t^{\alpha }$ (these are denoted by $x_{\alpha }^{i}$) and with respect to $%
s^{a}$ (these are denoted by $y_{a}^{i}$), respectively. The coordinate
transformation rules on $J^{1}(TM,N)$ are given by the general
transformation laws [9]%
\begin{equation*}
\left\{ 
\begin{array}{l}
\overline{T}^{A}=\overline{T}^{A}\left( T^{B}\right) \medskip \\ 
\overline{x}^{i}=\overline{x}^{i}\left( x^{j}\right) \medskip \\ 
\overline{X}_{A}^{i}=\dfrac{\partial \overline{x}^{i}}{\partial x^{j}}\dfrac{%
\partial T^{B}}{\partial \overline{T}^{A}}X_{B}^{j}.%
\end{array}%
\right.
\end{equation*}

Consequently, direct local computations say us that the local coordinates on 
$J^{1}(TM,N)$ are $\left( t^{\alpha },s^{a},x^{i},x_{\alpha
}^{i},y_{a}^{i}\right) $ and that they transform by the rules%
\begin{equation}
\left\{ 
\begin{array}{l}
\overline{t}^{\alpha }=\overline{t}^{\alpha }\left( t^{\beta }\right)
\medskip \\ 
\overline{s}^{a}=\dfrac{\partial \overline{t}^{a}}{\partial t^{b}}%
s^{b}\medskip \\ 
\overline{x}^{i}=\overline{x}^{i}\left( x^{j}\right) \medskip \\ 
\overline{x}_{\alpha }^{i}=\dfrac{\partial \overline{x}^{i}}{\partial x^{j}}%
\dfrac{\partial t^{\beta }}{\partial \overline{t}^{\alpha }}x_{\beta }^{j}+%
\dfrac{\partial \overline{x}^{i}}{\partial x^{j}}\dfrac{\partial ^{2}t^{b}}{%
\partial \overline{t}^{\alpha }\partial \overline{t}^{c}}\dfrac{\partial 
\overline{t}^{c}}{\partial t^{d}}s^{d}y_{b}^{j}\medskip \\ 
\overline{y}_{a}^{i}=\dfrac{\partial \overline{x}^{i}}{\partial x^{j}}\dfrac{%
\partial t^{b}}{\partial \overline{t}^{a}}y_{b}^{j}.%
\end{array}%
\right.  \label{TJ}
\end{equation}

Firstly, let us suppose that%
\begin{equation*}
\left( \Gamma _{AB}^{C}\right) =\left( \text{ }^{1}\Gamma _{\alpha \beta
}^{\gamma },\text{ }^{2}\Gamma _{\alpha \beta }^{c},\text{ }^{3}\Gamma
_{a\beta }^{\gamma },\text{ }^{4}\Gamma _{a\beta }^{c},\text{ }^{5}\Gamma
_{\alpha b}^{\gamma },\text{ }^{6}\Gamma _{\alpha b}^{c},\text{ }^{7}\Gamma
_{ab}^{\gamma },\text{ }^{8}\Gamma _{ab}^{c}\right)
\end{equation*}%
are the \textit{normal components} of a linear connection $\nabla $ on $%
TM\backslash \{0\}$, in the sense that we have%
\begin{equation*}
\nabla _{\dfrac{\partial }{\partial T^{B}}}\dfrac{\partial }{\partial T^{A}}%
=\Gamma _{AB}^{C}\dfrac{\partial }{\partial T^{C}},
\end{equation*}%
that is%
\begin{equation*}
\begin{array}{ll}
\nabla _{\dfrac{\partial }{\partial t^{\beta }}}\dfrac{\partial }{\partial
t^{\alpha }}=\text{ }^{1}\Gamma _{\alpha \beta }^{\gamma }\dfrac{\partial }{%
\partial t^{\gamma }}+\text{ }^{2}\Gamma _{\alpha \beta }^{c}\dfrac{\partial 
}{\partial s^{c}}, & \nabla _{\dfrac{\partial }{\partial t^{\beta }}}\dfrac{%
\partial }{\partial s^{a}}=\text{ }^{3}\Gamma _{a\beta }^{\gamma }\dfrac{%
\partial }{\partial t^{\gamma }}+\text{ }^{4}\Gamma _{a\beta }^{c}\dfrac{%
\partial }{\partial s^{c}},\medskip \\ 
\nabla _{\dfrac{\partial }{\partial s^{b}}}\dfrac{\partial }{\partial
t^{\alpha }}=\text{ }^{5}\Gamma _{\alpha b}^{\gamma }\dfrac{\partial }{%
\partial t^{\gamma }}+\text{ }^{6}\Gamma _{\alpha b}^{c}\dfrac{\partial }{%
\partial s^{c}}, & \nabla _{\dfrac{\partial }{\partial s^{b}}}\dfrac{%
\partial }{\partial s^{a}}=\text{ }^{7}\Gamma _{ab}^{\gamma }\dfrac{\partial 
}{\partial t^{\gamma }}+\text{ }^{8}\Gamma _{ab}^{c}\dfrac{\partial }{%
\partial s^{c}}.%
\end{array}%
\end{equation*}%
Then, following the geometrical ideas from [9, pp. 25], we underline that
the components%
\begin{equation}
M_{(B)A}^{(j)}=-\Gamma _{AB}^{C}X_{C}^{j},  \label{2stel}
\end{equation}%
where 
\begin{equation*}
X_{C}^{j}=\left( x_{\gamma }^{j},y_{c}^{j}\right) ,
\end{equation*}%
represent a \textit{temporal nonlinear connection} on the 1-jet space $%
J^{1}(TM,N)$, \textit{naturally attached to the linear connection} $\nabla $ 
\textit{from} $TM\backslash \{0\}.$

In this geometrical context, we point out that the Berwald adapted
components $B_{\alpha \beta }^{\gamma }=B_{\alpha \beta }^{\gamma }(t^{\mu
},s^{a})$ define a linear d-connection $^{\mathbf{b}}\nabla $ on $%
TM\backslash \{0\}$, given in the adapted basis (\ref{AB}) by the relations%
\begin{equation}
\begin{array}{ll}
^{\mathbf{b}}\nabla _{\dfrac{\delta }{\delta t^{\beta }}}\dfrac{\delta }{%
\delta t^{\alpha }}=B_{\alpha \beta }^{\gamma }\dfrac{\delta }{\delta
t^{\gamma }}, & ^{\mathbf{b}}\nabla _{\dfrac{\delta }{\delta t^{\beta }}}%
\dfrac{\partial }{\partial s^{a}}=B_{a\beta }^{c}\dfrac{\partial }{\partial
s^{c}},\medskip \\ 
^{\mathbf{b}}\nabla _{\dfrac{\partial }{\partial s^{b}}}\dfrac{\delta }{%
\delta t^{\alpha }}=0, & ^{\mathbf{b}}\nabla _{\dfrac{\partial }{\partial
s^{b}}}\dfrac{\partial }{\partial s^{a}}=0.%
\end{array}
\label{BA}
\end{equation}

\begin{remark}
The adapted torsion and curvature d-tensors of the Berwald linear
d-connection $^{\mathbf{b}}\nabla $ on $TM\backslash \{0\}$ are given by the
formulas (\ref{BT}) and (\ref{BC}).
\end{remark}

Taking into account the relations (\ref{BA}), together with the equality%
\begin{equation*}
\frac{\partial }{\partial t^{\alpha }}=\frac{\delta }{\delta t^{\alpha }}%
+N_{\alpha }^{b}\frac{\partial }{\partial s^{b}},
\end{equation*}%
by direct computations, we deduce that the Berwald d-connection $^{\mathbf{b}%
}\nabla $ on $TM\backslash \{0\}$ has, in the \textit{normal basis} 
\begin{equation*}
\left\{ \frac{\partial }{\partial t^{\alpha }},\frac{\partial }{\partial
s^{a}}\right\} \subset \mathcal{X}(TM\backslash \{0\}),
\end{equation*}%
the \textit{normal components} $\left( \text{ }^{\mathbf{b}}\Gamma
_{AB}^{C}\right) =$%
\begin{equation}
\begin{array}{l}
\text{ }^{1\mathbf{b}}\Gamma _{\alpha \beta }^{\gamma }=B_{\alpha \beta
}^{\gamma },\text{ }^{2\mathbf{b}}\Gamma _{\alpha \beta }^{c}=N_{\alpha
:\beta }^{c},\text{ }^{3\mathbf{b}}\Gamma _{a\beta }^{\gamma }=0,\text{ }^{4%
\mathbf{b}}\Gamma _{a\beta }^{c}=B_{a\beta }^{c},\medskip \\ 
\text{ }^{5\mathbf{b}}\Gamma _{\alpha b}^{\gamma }=0,\text{ }^{6\mathbf{b}%
}\Gamma _{\alpha b}^{c}=B_{\alpha b}^{c},\text{ }^{7\mathbf{b}}\Gamma
_{ab}^{\gamma }=0,\text{ }^{8\mathbf{b}}\Gamma _{ab}^{c}=0,%
\end{array}
\label{R1}
\end{equation}%
where%
\begin{equation*}
N_{\alpha :\beta }^{c}=\frac{\partial N_{\alpha }^{c}}{\partial t^{\beta }}%
+N_{\alpha }^{d}B_{d\beta }^{c}-N_{\gamma }^{c}B_{\alpha \beta }^{\gamma }.
\end{equation*}

As a consequence, via the formula (\ref{2stel}) applied to the Berwald
linear d-connection $^{\mathbf{b}}\nabla $ from $TM\backslash \{0\}$, we
find the following geometrical result and concept:

\begin{definition}
The set of local functions%
\begin{equation*}
\left( \text{ }^{\mathbf{b}}M_{(B)A}^{(j)}\right) =\left( \text{ }^{1\mathbf{%
b}}M_{(\beta )\alpha }^{(j)},\text{ }^{2\mathbf{b}}M_{(b)\alpha }^{(j)},%
\text{ }^{3\mathbf{b}}M_{(\beta )a}^{(j)},\text{ }^{4\mathbf{b}%
}M_{(b)a}^{(j)}\right) ,
\end{equation*}%
where%
\begin{equation}
\begin{array}{ll}
\text{ }^{1\mathbf{b}}M_{(\beta )\alpha }^{(j)}=-B_{\alpha \beta }^{\gamma
}x_{\gamma }^{j}-N_{\alpha :\beta }^{c}y_{c}^{j}, & \text{ }^{2\mathbf{b}%
}M_{(b)\alpha }^{(j)}=-B_{\alpha b}^{c}y_{c}^{j}\medskip \\ 
\text{ }^{3\mathbf{b}}M_{(\beta )a}^{(j)}=-B_{a\beta }^{c}y_{c}^{j}, & \text{
}^{4\mathbf{b}}M_{(b)a}^{(j)}=0,%
\end{array}
\label{R2}
\end{equation}%
represents a temporal nonlinear connection on $J^{1}(TM,N)$, which may be
called the \textbf{Berwald temporal nonlinear connection} on $J^{1}(TM,N)$.
\end{definition}

Secondly, let us consider that $\varphi :(M,F)\rightarrow (N,\widetilde{F})$
is an affine map between the Finsler manifolds $(M,F)$ and $(N,\widetilde{F}%
) $. Then, it is important to note that, under a change of coordinates on $M$
and $N$, the behaviour on $J^{1}(TM,N)$ of the coordinates $\left(
y_{a}^{i}\right) $ (see the relations (\ref{TJ})) is the same with that of
the components $\left( \varphi _{\alpha }^{i}\right) $. Consequently, by a
convenient jet extension, the Berwald adapted components 
\begin{equation*}
\widetilde{B}_{jk}^{i}=\widetilde{B}_{jk}^{i}(\varphi ^{l},\varphi
_{\varepsilon }^{l}s^{\varepsilon }),
\end{equation*}%
which appear in the equations of the affine maps (\ref{Af}), can be well
represented on the 1-jet space $J^{1}(TM,N)$\ by the geometrical objects 
\begin{equation*}
\widetilde{B}_{jk}^{i}=\widetilde{B}_{jk}^{i}(x^{l},y_{a}^{l}s^{a}),
\end{equation*}%
whose transformation rules on $J^{1}(TM,N)$ are (for more details, please
consult [6, pp. 122] or [3, pp. 43])%
\begin{equation}
\overline{\widetilde{B}}_{rs}^{q}=\widetilde{B}_{jk}^{i}\frac{\partial x^{j}%
}{\partial \overline{x}^{r}}\frac{\partial x^{k}}{\partial \overline{x}^{s}}%
\frac{\partial \overline{x}^{q}}{\partial x^{i}}+\frac{\partial ^{2}x^{i}}{%
\partial \overline{x}^{r}\partial \overline{x}^{s}}\frac{\partial \overline{x%
}^{q}}{\partial x^{i}}.  \label{Ber-Christ}
\end{equation}

Following the jet geometrical ideas from [9] or [10], it immediately follows
that the components [9, pp. 25]%
\begin{equation}
\text{ }^{\mathbf{b}}N_{(B)i}^{(j)}=\widetilde{B}_{ik}^{j}X_{B}^{k},
\label{1stea}
\end{equation}%
where 
\begin{equation*}
X_{B}^{k}=\left( x_{\beta }^{k},y_{b}^{k}\right) ,
\end{equation*}%
represent a \textit{spatial nonlinear connection} on the 1-jet space $%
J^{1}(TM,N)$. In conclusion, putting $B=\beta $ and $B=b$, respectively,
into the formula (\ref{1stea}), we can enunciate the following geometrical
result and concept:

\begin{definition}
The set of local functions%
\begin{equation*}
\left( \text{ }^{\mathbf{b}}N_{(B)i}^{(j)}\right) =\left( \text{ }^{1\mathbf{%
b}}N_{(\beta )i}^{(j)},\text{ }^{2\mathbf{b}}N_{(b)i}^{(j)}\right) ,
\end{equation*}%
where%
\begin{equation}
\begin{array}{ll}
\text{ }^{1\mathbf{b}}N_{(\beta )i}^{(j)}=\widetilde{B}_{ik}^{j}x_{\beta
}^{k}, & \text{ }^{2\mathbf{b}}N_{(b)i}^{(j)}=\widetilde{B}%
_{ik}^{j}y_{b}^{k},%
\end{array}
\label{R3}
\end{equation}%
represents a spatial nonlinear connection on $J^{1}(TM,N)$, which may be
called the \textbf{Berwald spatial nonlinear connection} on $J^{1}(TM,N)$.
\end{definition}

\begin{remark}
The set of local functions%
\begin{equation*}
\Gamma _{\text{jet}}^{\mathbf{b}}=\left( \text{ }^{\mathbf{b}}M_{(B)A}^{(j)},%
\text{ }^{\mathbf{b}}N_{(B)i}^{(j)}\right)
\end{equation*}%
is a nonlinear connection on $J^{1}(TM,N)$, which may be called the \textbf{%
jet Berwald nonlinear connection} on $J^{1}(TM,N)$.
\end{remark}

The Berwald nonlinear connection $\Gamma _{\text{\textit{jet}}}^{\mathbf{b}}$%
, whose local components are given by (\ref{R2}) and (\ref{R3}), produces
the \textit{jet adapted basis} [9, pp. 24]%
\begin{equation*}
\left\{ \frac{\delta ^{\mathbf{J}}}{\delta T^{A}},\frac{\delta ^{\mathbf{J}}%
}{\delta x^{i}},\frac{\partial }{\partial X_{A}^{i}}\right\} \subset 
\mathcal{X}(J^{1}(TM,N)),
\end{equation*}%
where%
\begin{equation*}
\begin{array}{cc}
\dfrac{\delta ^{\mathbf{J}}}{\delta T^{A}}=\dfrac{\partial }{\partial T^{A}}-%
\text{ }^{\mathbf{b}}M_{(B)A}^{(j)}\dfrac{\partial }{\partial X_{B}^{j}}, & 
\dfrac{\delta ^{\mathbf{J}}}{\delta x^{i}}=\dfrac{\partial }{\partial x^{i}}-%
\text{ }^{\mathbf{b}}N_{(B)i}^{(j)}\dfrac{\partial }{\partial X_{B}^{j}}.%
\end{array}%
\end{equation*}

Taking into account that $A=(\alpha ,a)$ and using the formulas (\ref{R2})
and (\ref{R3}), it is easy to deduce

\begin{proposition}
The elements of the \textbf{jet adapted basis} are%
\begin{equation*}
\left\{ \frac{\delta ^{\mathbf{J}}}{\delta t^{\alpha }},\frac{\delta ^{%
\mathbf{J}}}{\delta s^{a}},\frac{\delta ^{\mathbf{J}}}{\delta x^{i}},\frac{%
\partial }{\partial x_{\alpha }^{i}},\frac{\partial }{\partial y_{a}^{i}}%
\right\} \subset \mathcal{X}(J^{1}(TM,N)),
\end{equation*}%
where%
\begin{equation}
\begin{array}{l}
\dfrac{\delta ^{\mathbf{J}}}{\delta t^{\alpha }}=\dfrac{\partial }{\partial
t^{\alpha }}+\left( B_{\alpha \beta }^{\gamma }x_{\gamma }^{j}+N_{\alpha
:\beta }^{c}y_{c}^{j}\right) \dfrac{\partial }{\partial x_{\beta }^{j}}%
+B_{\alpha b}^{c}y_{c}^{j}\dfrac{\partial }{\partial y_{b}^{j}},\medskip \\ 
\dfrac{\delta ^{\mathbf{J}}}{\delta s^{a}}=\dfrac{\partial }{\partial s^{a}}%
+B_{a\beta }^{c}y_{c}^{j}\dfrac{\partial }{\partial x_{\beta }^{j}},\medskip
\\ 
\dfrac{\delta ^{\mathbf{J}}}{\delta x^{i}}=\dfrac{\partial }{\partial x^{i}}-%
\widetilde{B}_{ik}^{j}x_{\beta }^{k}\dfrac{\partial }{\partial x_{\beta }^{j}%
}-\widetilde{B}_{ik}^{j}y_{b}^{k}\dfrac{\partial }{\partial y_{b}^{j}}.%
\end{array}
\label{Der}
\end{equation}
\end{proposition}

Following again the jet geometrical ideas exposed in [9] and [10], the
Berwald linear d-connection $^{\mathbf{b}}\nabla $ from $TM\backslash \{0\}$%
, together with the Berwald linear d-connection $^{\mathbf{b}}\widetilde{%
\nabla }$ from $TN\backslash \{0\}$, produce a \textit{jet Berwald linear
d-connection }$B\Gamma _{\text{\textit{jet}}}^{\mathbf{b}}$ on $J^{1}(TM,N)$%
, taking as its \textit{jet adapted components} the following coefficients
[9, pp. 30]:

\begin{equation*}
\begin{array}{l}
B\Gamma _{\text{\textit{jet}}}^{\mathbf{b}}=\left( \text{ }\overline{G}%
_{BC}^{A}=\text{ }^{\mathbf{b}}\Gamma _{BC}^{A},\text{ }G_{iC}^{k}=0,\text{ }%
G_{(A)(j)C}^{(i)(B)}=-\delta _{j}^{i}\text{ }^{\mathbf{b}}\Gamma
_{CA}^{B},\right. \medskip \\ 
\text{\ \ \ \ \ \ \ \ \ \ \ \ }\left. \text{ }\overline{L}_{Bj}^{A}=0,\text{ 
}L_{ij}^{k}=\widetilde{B}_{ij}^{k},\text{ }L_{(A)(j)k}^{(i)(B)}=\delta
_{A}^{B}\widetilde{B}_{jk}^{i},\text{ }\right. \medskip \\ 
\text{\ \ \ \ \ \ \ \ \ \ \ \ }\left. \text{ }\overline{C}_{B(k)}^{A(C)}=0,%
\text{ }C_{i(k)}^{j(C)}=0,\text{ }C_{(A)(j)(k)}^{(i)(B)(C)}=0\text{ }\right)
,%
\end{array}%
\end{equation*}%
where%
\begin{equation*}
\delta _{A}^{B}=\left\{ 
\begin{array}{ll}
\delta _{\alpha }^{\beta }, & \text{if }A=\alpha ,\text{ }B=\beta \medskip
\\ 
\delta _{a}^{b}, & \text{if }A=a,\text{ }B=b\medskip \\ 
0, & \text{otherwise.}%
\end{array}%
\right.
\end{equation*}

\begin{remark}
The jet \textit{Berwald linear d-connection }$B\Gamma _{\text{\textit{jet}}%
}^{\mathbf{b}}$ is a $\Gamma _{\text{\textit{jet}}}^{\mathbf{b}}$-linear
connection on $J^{1}(TM,N)$. For more details, please consult [9, pp. 28] or
[10].
\end{remark}

Consequently, using the relations (\ref{R1}) and taking into account that
the indices $A,B,C,...$ have the form $(\alpha ,a),(\beta ,b),(\gamma
,c),... $, we obtain

\begin{proposition}
The nonvanishing jet adapted components of the \textbf{jet Berwald linear
d-connection} $B\Gamma _{\text{\textit{jet}}}^{\mathbf{b}}$ are only the
following \textbf{eleven} components:%
\begin{equation}
\begin{array}{l}
B\Gamma _{\text{\textit{jet}}}^{\mathbf{b}}=\left( \text{ }\overline{G}%
_{\beta \gamma }^{\alpha }=B_{\beta \gamma }^{\alpha },\text{ }\overline{G}%
_{\beta \gamma }^{a}=N_{\beta :\gamma }^{a},\text{ }\overline{G}_{b\gamma
}^{a}=B_{b\gamma }^{a},\text{ }\overline{G}_{\beta c}^{a}=B_{\beta
c}^{a},\right. \medskip \\ 
\text{ \ \ \ \ \ \ \ \ \ \ }\left. \text{ }G_{(\alpha )(j)\gamma
}^{(i)(\beta )}=-\delta _{j}^{i}B_{\gamma \alpha }^{\beta },\text{ }%
G_{(\alpha )(j)\gamma }^{(i)(b)}=-\delta _{j}^{i}N_{\gamma :\alpha
}^{b},\right. \medskip \\ 
\text{ \ \ \ \ \ \ \ \ \ \ }\left. \text{ }G_{(\alpha
)(j)c}^{(i)(b)}=-\delta _{j}^{i}B_{c\alpha }^{b},\text{ }G_{(a)(j)\gamma
}^{(i)(b)}=-\delta _{j}^{i}B_{\gamma a}^{b},\text{ }L_{ij}^{k}=\widetilde{B}%
_{ij}^{k},\text{ }\right. \medskip \\ 
\text{ \ \ \ \ \ \ \ \ \ \ }\left. \text{ }L_{(\alpha )(j)k}^{(i)(\beta
)}=\delta _{\alpha }^{\beta }\widetilde{B}_{jk}^{i},\text{ }%
L_{(a)(j)k}^{(i)(b)}=\delta _{a}^{b}\widetilde{B}_{jk}^{i}\text{ }\right) .%
\end{array}
\label{R1jet}
\end{equation}
\end{proposition}

Because the Riemann-Lagrange geometry of the general $\Gamma $-linear
connections on 1-jet spaces, in the sense of their d-torsions and
d-curvatures, is now completely done in [9] and [10], it follows that we can
compute on $J^{1}(TM,N)$ the adapted components of the torsion and curvature
d-tensors produced by the Berwald linear d-connection $B\Gamma _{\text{%
\textit{jet}}}^{\mathbf{b}}$. In a such jet Riemann-Lagrange geometrical
context, using the formulas (\ref{BT}) and (\ref{BC}), we can give the
following geometrical results:

\begin{theorem}
The jet Berwald linear d-connection $B\Gamma _{\text{\textit{jet}}}^{\mathbf{%
b}}$ on $J^{1}(TM,N)$ is characterized by \textbf{fifteen }nonvanishing
local adapted \textbf{d-torsions}:\medskip

\textbf{(T1)} $T_{\beta \gamma }^{a}=N_{\beta :\gamma }^{a}-N_{\gamma :\beta
}^{a}=$ $^{\mathbf{b}}\mathcal{R}_{\beta \gamma }^{\alpha }$,\medskip

\textbf{(T2)} $P_{(\mu )i(j)}^{(m)\text{ }(b)}=$ $^{\mathbf{b}}\widetilde{%
\mathcal{P}}_{ikj}^{m}x_{\mu }^{k}s^{b}$,\medskip

\textbf{(T3)} $P_{(c)i(j)}^{(m)\text{ }(b)}=$ $^{\mathbf{b}}\widetilde{%
\mathcal{P}}_{ikj}^{m}y_{c}^{k}s^{b}$,\medskip

\textbf{(T4)} $R_{(\mu )\alpha \beta }^{(m)}=-\left[ ^{\mathbf{b}}\mathcal{R}%
_{\mu \alpha \beta }^{\varepsilon }+\text{ }^{\mathbf{b}}\mathcal{P}_{\mu
\alpha c}^{\varepsilon }N_{\beta }^{c}-\text{ }^{\mathbf{b}}\mathcal{P}_{\mu
\beta c}^{\varepsilon }N_{\alpha }^{c}\right] x_{\varepsilon }^{m}+\medskip$

$\ \ \ \ \ \ \ $\ $\ \ \ \ \ \ \ \ \ \ \ +\mathcal{A}_{\{\alpha ,\beta \}}%
\left[ \dfrac{\partial N_{\beta :\mu }^{c}}{\partial t^{\alpha }}+N_{\beta
:\mu }^{d}B_{d\alpha }^{c}-N_{\beta :\gamma }^{c}B_{\alpha \mu }^{\gamma }%
\right] y_{c}^{m}$,\medskip

\textbf{(T5)} $R_{(\mu )\alpha b}^{(m)}=-$ $^{\mathbf{b}}\mathcal{P}_{\mu
\alpha b}^{\varepsilon }x_{\varepsilon }^{m}+\left[ \dfrac{\partial B_{b\mu
}^{c}}{\partial t^{\alpha }}-\dfrac{\partial N_{\alpha :\mu }^{c}}{\partial
s^{b}}+B_{b\mu }^{d}B_{d\alpha }^{c}-B_{\alpha \mu }^{\gamma }B_{\gamma
b}^{c}\right] y_{c}^{m}$,\medskip

\textbf{(T6)} $R_{(\mu )a\beta }^{(m)}=$ $^{\mathbf{b}}\mathcal{P}_{\mu
a\beta }^{\varepsilon }x_{\varepsilon }^{m}-\left[ \dfrac{\partial B_{a\mu
}^{c}}{\partial t^{\beta }}-\dfrac{\partial N_{\beta :\mu }^{c}}{\partial
s^{a}}+B_{a\mu }^{d}B_{d\beta }^{c}-B_{\beta \mu }^{\gamma }B_{\gamma a}^{c}%
\right] y_{c}^{m}$,\medskip

\textbf{(T7)} $R_{(c)\alpha \beta }^{(m)}=-\left[ ^{\mathbf{b}}\mathcal{R}%
_{c\alpha \beta }^{d}+\text{ }^{\mathbf{b}}\mathcal{P}_{c\alpha
f}^{d}N_{\beta }^{f}-\text{ }^{\mathbf{b}}\mathcal{P}_{c\beta
f}^{d}N_{\alpha }^{f}\right] y_{d}^{m}$,$\medskip$

\textbf{(T8)} $R_{(c)\alpha b}^{(m)}=-$ $^{\mathbf{b}}\mathcal{P}_{c\alpha
b}^{d}y_{d}^{m}$,\medskip

\textbf{(T9)} $R_{(c)a\beta }^{(m)}=$ $^{\mathbf{b}}\mathcal{P}_{ca\beta
}^{d}y_{d}^{m}$,\medskip

\textbf{(T10)} $R_{(\mu )\alpha j}^{(m)}=-$ $^{\mathbf{b}}\widetilde{%
\mathcal{P}}_{jkl}^{m}B_{\alpha b}^{c}s^{b}x_{\mu }^{k}y_{c}^{l}$,\medskip

\textbf{(T11)} $R_{(\mu )aj}^{(m)}=-$ $^{\mathbf{b}}\widetilde{\mathcal{P}}%
_{jkl}^{m}x_{\mu }^{k}y_{a}^{l}$,\medskip\ 

\textbf{(T12)} $R_{(c)\alpha j}^{(m)}=-$ $^{\mathbf{b}}\widetilde{\mathcal{P}%
}_{jkl}^{m}B_{\alpha b}^{d}s^{b}y_{c}^{k}y_{d}^{l}$,\medskip

\textbf{(T13)} $R_{(c)aj}^{(m)}=-$ $^{\mathbf{b}}\widetilde{\mathcal{P}}%
_{jkl}^{m}y_{c}^{k}y_{a}^{l}$,\medskip

\textbf{(T14)} $R_{(\mu )ij}^{(m)}=\left[ \text{ }^{\mathbf{b}}\widetilde{%
\mathcal{R}}_{kij}^{m}+\text{ }^{\mathbf{b}}\widetilde{\mathcal{P}}_{kil}^{m}%
\widetilde{N}_{j}^{l}-\text{ }^{\mathbf{b}}\widetilde{\mathcal{P}}_{kjl}^{m}%
\widetilde{N}_{i}^{l}\right] x_{\mu }^{k}-\medskip$

$\ \ \ \ \ \ \ \ $\ $\ \ \ \ \ \ \ \ \ -\left[ \text{ }^{\mathbf{b}}%
\widetilde{\mathcal{P}}_{kil}^{m}\widetilde{B}_{jp}^{l}-\text{ }^{\mathbf{b}}%
\widetilde{\mathcal{P}}_{kjl}^{m}\widetilde{B}_{ip}^{l}\right] x_{\mu
}^{k}y_{a}^{p}s^{a}$,\medskip

\textbf{(T15)} $R_{(c)ij}^{(m)}=\left[ \text{ }^{\mathbf{b}}\widetilde{%
\mathcal{R}}_{kij}^{m}+\text{ }^{\mathbf{b}}\widetilde{\mathcal{P}}_{kil}^{m}%
\widetilde{N}_{j}^{l}-\text{ }^{\mathbf{b}}\widetilde{\mathcal{P}}_{kjl}^{m}%
\widetilde{N}_{i}^{l}\right] y_{c}^{k}-\medskip$

$\ \ \ \ \ \ \ \ $\ $\ \ \ \ \ \ \ \ \ -\left[ \text{ }^{\mathbf{b}}%
\widetilde{\mathcal{P}}_{kil}^{m}\widetilde{B}_{jp}^{l}-\text{ }^{\mathbf{b}}%
\widetilde{\mathcal{P}}_{kjl}^{m}\widetilde{B}_{ip}^{l}\right]
y_{c}^{k}y_{a}^{p}s^{a}$,\medskip

where $\mathcal{A}_{\{\alpha ,\beta \}}$ means an alternate sum;
\end{theorem}

\begin{proof}
The nonvanishing local adapted d-torsions of the Berwald linear d-connection 
$B\Gamma _{\text{\textit{jet}}}^{\mathbf{b}}$ on $J^{1}(TM,N)$ are given by
the general formulas [9, pp. 34]:\medskip

\textbf{(t1)} $T_{AB}^{M}=\overline{G}_{AB}^{M}-\overline{G}%
_{BA}^{M},\medskip$

\textbf{(t2) }$P_{(M)A(j)}^{(m)\text{ \ }(B)}=\dfrac{\partial \left[ \text{ }%
^{\mathbf{b}}M_{(M)A}^{(m)}\right] }{\partial X_{B}^{j}}%
-G_{(M)(j)A}^{(m)(B)},\medskip$

\textbf{(t3) }$P_{(M)i(j)}^{(m)\text{ }(B)}=\dfrac{\partial \left[ \text{ }^{%
\mathbf{b}}N_{(M)i}^{(m)}\right] }{\partial X_{B}^{j}}-L_{(M)(j)i}^{(m)(B)},%
\medskip$

\textbf{(t4) }$R_{(M)AB}^{(m)}=\dfrac{\delta ^{\mathbf{J}}\left[ \text{ }^{%
\mathbf{b}}M_{(M)A}^{(m)}\right] }{\delta T^{B}}-\dfrac{\delta ^{\mathbf{J}}%
\left[ \text{ }^{\mathbf{b}}M_{(M)B}^{(m)}\right] }{\delta T^{A}},\medskip$

\textbf{(t5) }$R_{(M)Aj}^{(m)}=\dfrac{\delta ^{\mathbf{J}}\left[ \text{ }^{%
\mathbf{b}}M_{(M)A}^{(m)}\right] }{\delta x^{j}}-\dfrac{\delta ^{\mathbf{J}}%
\left[ \text{ }^{\mathbf{b}}N_{(M)j}^{(m)}\right] }{\delta T^{A}},\medskip$

\textbf{(t6) }$R_{(M)ij}^{(m)}=\dfrac{\delta ^{\mathbf{J}}\left[ \text{ }^{%
\mathbf{b}}N_{(M)i}^{(m)}\right] }{\delta x^{j}}-\dfrac{\delta ^{\mathbf{J}}%
\left[ \text{ }^{\mathbf{b}}N_{(M)j}^{(m)}\right] }{\delta x^{i}},\medskip$

where $\left( T^{A}\right) =\left( t^{\alpha },s^{a}\right) $ and $\left(
X_{A}^{i}\right) =\left( x_{\alpha }^{i},y_{a}^{i}\right) $.

Taking into account that the indices $A,B,...$ are indices of kind $(\alpha
,a),(\beta ,b),...$ and using the fomulas (\ref{R1jet}), (\ref{R2}), (\ref%
{R3}) and (\ref{Der}), by laborious local computations, we find the required
result.
\end{proof}

\begin{theorem}
The jet Berwald linear d-connection $B\Gamma _{\text{\textit{jet}}}^{\mathbf{%
b}}$ on $J^{1}(TM,N)$ is characterized by \textbf{thirty }nonvanishing local
adapted \textbf{d-curvatures}:\medskip

\textbf{(C1)} $\overline{R}_{\alpha \beta \gamma }^{\delta }=$ $^{\mathbf{b}}%
\mathcal{R}_{\alpha \beta \gamma }^{\delta }+$ $^{\mathbf{b}}\mathcal{P}%
_{\alpha \beta c}^{\delta }N_{\gamma }^{c}-$ $^{\mathbf{b}}\mathcal{P}%
_{\alpha \gamma c}^{\delta }N_{\beta }^{c}$,\medskip

\textbf{(C2) }$\overline{R}_{\alpha \beta \gamma }^{d}=\mathcal{A}_{\{\beta
,\gamma \}}\left[ \dfrac{\partial N_{\alpha :\beta }^{d}}{\partial t^{\gamma
}}+N_{\alpha :\beta }^{c}B_{c\gamma }^{d}-N_{\mu :\beta }^{d}B_{\alpha
\gamma }^{\mu }\right] $,\medskip

\textbf{(C3)} $\overline{R}_{\alpha \beta c}^{\delta }=$ $^{\mathbf{b}}%
\mathcal{P}_{\alpha \beta c}^{\delta }$,\medskip

\textbf{(C4) }$\overline{R}_{\alpha \beta c}^{d}=\dfrac{\partial N_{\alpha
:\beta }^{d}}{\partial s^{c}}-\dfrac{\partial B_{\alpha c}^{d}}{\partial
t^{\beta }}+B_{\alpha \beta }^{\mu }B_{\mu c}^{d}-B_{\alpha c}^{f}B_{f\beta
}^{d}$,\medskip

\textbf{(C5)} $\overline{R}_{\alpha b\gamma }^{\delta }=-$ $^{\mathbf{b}}%
\mathcal{P}_{\alpha b\gamma }^{\delta }$,\medskip

\textbf{(C6)} $\overline{R}_{\alpha b\gamma }^{d}=\dfrac{\partial B_{\alpha
b}^{d}}{\partial t^{\gamma }}-\dfrac{\partial N_{\alpha :\gamma }^{d}}{%
\partial s^{b}}+B_{\alpha b}^{c}B_{c\gamma }^{d}-B_{\alpha \gamma }^{\mu
}B_{\mu b}^{d}$,\medskip

\textbf{(C7)} $\overline{R}_{a\beta \gamma }^{d}=$ $^{\mathbf{b}}\mathcal{R}%
_{a\beta \gamma }^{d}+$ $^{\mathbf{b}}\mathcal{P}_{a\beta c}^{d}N_{\gamma
}^{c}-$ $^{\mathbf{b}}\mathcal{P}_{a\gamma c}^{d}N_{\beta }^{c}$,\medskip

\textbf{(C8)} $\overline{R}_{a\beta c}^{d}=$ $^{\mathbf{b}}\mathcal{P}%
_{a\beta c}^{d}$,\medskip

\textbf{(C9)} $\overline{R}_{ab\gamma }^{d}=-$ $^{\mathbf{b}}\mathcal{P}%
_{ab\gamma }^{d}$,\medskip

\textbf{(C10)} $R_{i\beta k}^{l}=-$ $^{\mathbf{b}}\widetilde{\mathcal{P}}%
_{ikj}^{l}B_{\beta a}^{c}s^{a}y_{c}^{j}$,\medskip

\textbf{(C11)} $R_{ibk}^{l}=-$ $^{\mathbf{b}}\widetilde{\mathcal{P}}%
_{ikj}^{l}y_{b}^{j}$,\medskip

\textbf{(C12)} $R_{ijk}^{l}=$ $^{\mathbf{b}}\widetilde{\mathcal{R}}%
_{ijk}^{l}+$ $^{\mathbf{b}}\widetilde{\mathcal{P}}_{ijr}^{l}\widetilde{N}%
_{k}^{r}-$ $^{\mathbf{b}}\widetilde{\mathcal{P}}_{ikr}^{l}\widetilde{N}%
_{j}^{r}+\left[ \text{ }^{\mathbf{b}}\widetilde{\mathcal{P}}_{ikr}^{l}%
\widetilde{B}_{jp}^{r}-\text{ }^{\mathbf{b}}\widetilde{\mathcal{P}}_{ijr}^{l}%
\widetilde{B}_{kp}^{r}\right] y_{b}^{p}s^{b}$,$\medskip $

\textbf{(C13)} $P_{ij(k)}^{l\text{ }(c)}=$ $^{\mathbf{b}}\widetilde{\mathcal{%
P}}_{ijk}^{l}s^{c}$,\medskip

\textbf{(C14) }$R_{(\varepsilon )(i)\beta \gamma }^{(l)(\alpha )}=-\delta
_{i}^{l}\cdot \left[ \text{ }^{\mathbf{b}}\mathcal{R}_{\varepsilon \beta
\gamma }^{\alpha }+\text{ }^{\mathbf{b}}\mathcal{P}_{\varepsilon \beta
c}^{\alpha }N_{\gamma }^{c}-\text{ }^{\mathbf{b}}\mathcal{P}_{\varepsilon
\gamma c}^{\alpha }N_{\beta }^{c}\right] =-\delta _{i}^{l}\cdot \overline{R}%
_{\varepsilon \beta \gamma }^{\alpha }$,\medskip

\textbf{(C15) }$R_{(\varepsilon )(i)\beta c}^{(l)(\alpha )}=-\delta
_{i}^{l}\cdot $ $^{\mathbf{b}}\mathcal{P}_{\varepsilon \beta c}^{\alpha
}=-\delta _{i}^{l}\cdot \overline{R}_{\varepsilon \beta c}^{\alpha }$%
,\medskip

\textbf{(C16) }$R_{(\varepsilon )(i)b\gamma }^{(l)(\alpha )}=\delta
_{i}^{l}\cdot $ $^{\mathbf{b}}\mathcal{P}_{\varepsilon b\gamma }^{\alpha
}=-\delta _{i}^{l}\cdot \overline{R}_{\varepsilon b\gamma }^{\alpha }$%
,\medskip

\textbf{(C17) }$R_{(\varepsilon )(i)\beta \gamma }^{(l)(a)}=-\delta
_{i}^{l}\cdot \mathcal{A}_{\{\beta ,\gamma \}}\left[ \dfrac{\partial
N_{\beta :\varepsilon }^{a}}{\partial t^{\gamma }}+N_{\beta :\varepsilon
}^{c}B_{c\gamma }^{a}-N_{\beta :\mu }^{a}B_{\varepsilon \gamma }^{\mu }%
\right] $,\medskip

\textbf{(C18) }$R_{(\varepsilon )(i)\beta c}^{(l)(a)}=-\delta _{i}^{l}\cdot %
\left[ \dfrac{\partial N_{\beta :\varepsilon }^{a}}{\partial s^{c}}-\dfrac{%
\partial B_{c\varepsilon }^{a}}{\partial t^{\beta }}+B_{\beta \varepsilon
}^{\mu }B_{\mu c}^{a}-B_{c\varepsilon }^{d}B_{d\beta }^{a}\right] $,\medskip

\textbf{(C19) }$R_{(\varepsilon )(i)b\gamma }^{(l)(a)}=-\delta _{i}^{l}\cdot %
\left[ \dfrac{\partial B_{b\varepsilon }^{a}}{\partial t^{\gamma }}-\dfrac{%
\partial N_{\gamma :\varepsilon }^{a}}{\partial s^{b}}+B_{b\varepsilon
}^{c}B_{c\gamma }^{a}-B_{\gamma \varepsilon }^{\mu }B_{\mu b}^{a}\right] $%
,\medskip

\textbf{(C20) }$R_{(d)(i)\beta \gamma }^{(l)(a)}=-\delta _{i}^{l}\cdot \left[
\text{ }^{\mathbf{b}}\mathcal{R}_{d\beta \gamma }^{a}+\text{ }^{\mathbf{b}}%
\mathcal{P}_{d\beta c}^{a}N_{\gamma }^{c}-\text{ }^{\mathbf{b}}\mathcal{P}%
_{d\gamma c}^{a}N_{\beta }^{c}\right] =-\delta _{i}^{l}\cdot \overline{R}%
_{d\beta \gamma }^{a}$,\medskip

\textbf{(C21) }$R_{(d)(i)\beta c}^{(l)(a)}=-\delta _{i}^{l}\cdot $ $^{%
\mathbf{b}}\mathcal{P}_{d\beta c}^{a}=-\delta _{i}^{l}\cdot \overline{R}%
_{d\beta c}^{a}$,\medskip

\textbf{(C22) }$R_{(d)(i)b\gamma }^{(l)(a)}=\delta _{i}^{l}\cdot $ $^{%
\mathbf{b}}\mathcal{P}_{db\gamma }^{a}=-\delta _{i}^{l}\cdot \overline{R}%
_{db\gamma }^{a}$,\medskip

\textbf{(C23) }$R_{(\varepsilon )(i)\beta k}^{(l)(\alpha )}=-\delta
_{\varepsilon }^{\alpha }\cdot $ $^{\mathbf{b}}\widetilde{\mathcal{P}}%
_{ikj}^{l}B_{\beta b}^{c}s^{b}y_{c}^{j}=\delta _{\varepsilon }^{\alpha
}\cdot R_{i\beta k}^{l}$,\medskip

\textbf{(C24) }$R_{(\varepsilon )(i)bk}^{(l)(\alpha )}=-\delta _{\varepsilon
}^{\alpha }\cdot $ $^{\mathbf{b}}\widetilde{\mathcal{P}}_{ikj}^{l}y_{b}^{j}=%
\delta _{\varepsilon }^{\alpha }\cdot R_{ibk}^{l}$,\medskip

\textbf{(C25) }$R_{(d)(i)\beta k}^{(l)(a)}=-\delta _{d}^{a}\cdot $ $^{%
\mathbf{b}}\widetilde{\mathcal{P}}_{ikj}^{l}B_{\beta
b}^{c}s^{b}y_{c}^{j}=\delta _{d}^{a}\cdot R_{i\beta k}^{l}$,\medskip

\textbf{(C26) }$R_{(d)(i)bk}^{(l)(a)}=-\delta _{d}^{a}\cdot $ $^{\mathbf{b}}%
\widetilde{\mathcal{P}}_{ikj}^{l}y_{b}^{j}=\delta _{d}^{a}\cdot R_{ibk}^{l}$%
,\medskip

\textbf{(C27) }$R_{(\varepsilon )(i)jk}^{(l)(\alpha )}=\delta _{\varepsilon
}^{\alpha }\cdot R_{ijk}^{l}$,\medskip

\textbf{(C28) }$R_{(d)(i)jk}^{(l)(a)}=\delta _{d}^{a}\cdot R_{ijk}^{l}$%
,\medskip

\textbf{(C29)} $P_{(\varepsilon )(i)j(k)}^{(l)(\alpha )\text{ }(c)}=\delta
_{\varepsilon }^{\alpha }\cdot $ $^{\mathbf{b}}\widetilde{\mathcal{P}}%
_{ijk}^{l}s^{c}=\delta _{\varepsilon }^{\alpha }\cdot P_{ij(k)}^{l\text{ }%
(c)}$,\medskip

\textbf{(C30)} $P_{(d)(i)j(k)}^{(l)(a)\text{ }(c)}=\delta _{d}^{a}\cdot $ $^{%
\mathbf{b}}\widetilde{\mathcal{P}}_{ijk}^{l}s^{c}=\delta _{d}^{a}\cdot
P_{ij(k)}^{l\text{ }(c)}$.
\end{theorem}

\begin{proof}
The nonvanishing local adapted d-curvatures of the Berwald linear
d-connection $B\Gamma _{\text{\textit{jet}}}^{\mathbf{b}}$ on $J^{1}(TM,N)$
are given by the general formulas [9, pp. 36]:\medskip

\textbf{(c1)} $\overline{R}_{ABC}^{D}=\dfrac{\delta ^{\mathbf{J}}\overline{G}%
_{AB}^{D}}{\delta T^{C}}-\dfrac{\delta ^{\mathbf{J}}\overline{G}_{AC}^{D}}{%
\delta T^{B}}+\overline{G}_{AB}^{M}\overline{G}_{MC}^{D}-\overline{G}%
_{AC}^{M}\overline{G}_{MB}^{D},\medskip$

\textbf{(c2) }$R_{iBk}^{l}=-\dfrac{\delta ^{\mathbf{J}}L_{ik}^{l}}{\delta
T^{B}},\medskip$

\textbf{(c3)} $R_{ijk}^{l}=\dfrac{\delta ^{\mathbf{J}}L_{ij}^{l}}{\delta
x^{k}}-\dfrac{\delta ^{\mathbf{J}}L_{ik}^{l}}{\delta x^{j}}%
+L_{ij}^{m}L_{mk}^{l}-L_{ik}^{m}L_{mj}^{l},\medskip$

\textbf{(c4)} $P_{ij(k)}^{l\text{ }(G)}=\dfrac{\partial L_{ij}^{l}}{\partial
X_{G}^{k}},\medskip$

\textbf{(c5)} $R_{(D)(i)BC}^{(l)(A)}=\dfrac{\delta ^{\mathbf{J}%
}G_{(D)(i)B}^{(l)(A)}}{\delta T^{C}}-\dfrac{\delta ^{\mathbf{J}%
}G_{(D)(i)C}^{(l)(A)}}{\delta T^{B}}+\medskip$

$\ \ \ \ \ \ \ \ \ \ \ \ \ \ \ \ \ \ \ \ \ \ \ \ +G_{(M)(i)B}^{(m)(A)}\cdot
G_{(D)(m)C}^{(l)(M)}-G_{(M)(i)C}^{(m)(A)}\cdot G_{(D)(m)B}^{(l)(M)},\medskip$

\textbf{(c6)} $R_{(D)(i)Bk}^{(l)(A)}=\dfrac{\delta ^{\mathbf{J}%
}G_{(D)(i)B}^{(l)(A)}}{\delta x^{k}}-\dfrac{\delta ^{\mathbf{J}%
}L_{(D)(i)k}^{(l)(A)}}{\delta T^{B}}+\medskip$

$\ \ \ \ \ \ \ \ \ \ \ \ \ \ \ \ \ \ \ \ \ \ \ \ +G_{(M)(i)B}^{(m)(A)}\cdot
L_{(D)(m)k}^{(l)(M)}-L_{(M)(i)k}^{(m)(A)}\cdot G_{(D)(m)B}^{(l)(M)},\medskip$

\textbf{(c7)} $R_{(D)(i)jk}^{(l)(A)}=\dfrac{\delta ^{\mathbf{J}%
}L_{(D)(i)j}^{(l)(A)}}{\delta x^{k}}-\dfrac{\delta ^{\mathbf{J}%
}L_{(D)(i)k}^{(l)(A)}}{\delta x^{j}}+\medskip$

$\ \ \ \ \ \ \ \ \ \ \ \ \ \ \ \ \ \ \ \ \ \ \ \ +L_{(M)(i)j}^{(m)(A)}\cdot
L_{(D)(m)k}^{(l)(M)}-L_{(M)(i)k}^{(m)(A)}\cdot L_{(D)(m)j}^{(l)(M)},\medskip$

\textbf{(c8)} $P_{(D)(i)j(k)}^{(l)(A)\text{ }(G)}=\dfrac{\partial
L_{(D)(i)j}^{(l)(A)}}{\partial X_{G}^{k}},\medskip$

where $\left( T^{B}\right) =\left( t^{\beta },s^{b}\right) $ and $\left(
X_{G}^{k}\right) =\left( x_{\gamma }^{k},y_{c}^{k}\right) $.

Taking into account that the indices $A,B,...$ are indices of kind $(\alpha
,a),(\beta ,b),...$ and using the fomulas (\ref{R1jet}) and (\ref{Der}), by
laborious local computations, we find what we were looking for.
\end{proof}

\begin{remark}[Open Problem]
In order to obtain geometrical informations on $J^{1}(TM,N)$ about our
starting affine map $\varphi :(M,F)\rightarrow (N,\widetilde{F})$, we can
replace again $y_{\alpha }^{i}$ with $\varphi _{\alpha }^{i}$. In this jet
geometrical context, the nondegenerate affine map $\varphi $ is effectively 
\textbf{"characterized"} by \textbf{eight} \textbf{jet d-torsions} \textbf{(}%
we refer to \textbf{(T2)},\textbf{\ (T3)},\textbf{\ (T10)-(T15))} and 
\textbf{twelve} \textbf{jet d-curvatures} \textbf{(}we refer to \textbf{%
(C10)-(C13)},\textbf{\ (C23)-(C30))}. It is an open problem what is the real
geometrical meaning of this intimate connection between the affine maps
between two Finsler manifolds and their attached d-torsions and d-curvatures
on the 1-jet space $J^{1}(TM,N)$. This open problem which is in author's
attention.
\end{remark}

\textbf{Acknowledgements. }I would like to thank to Professor S. Ianu\c{s}
who suggested me to elaborate a such geometrical study and who offered me
bibliographical notes in this direction.

\textbf{Author's address: }Mircea NEAGU, Str. L\u{a}m\^{a}i\c{t}ei, Nr. 66,
Bl. 93, Sc. G, Ap. 10, Bra\c{s}ov, BV 500371, Romania

\textbf{E-mail}: mirceaneagu73@yahoo.com

\begin{center}
University "Transilvania" of Bra\c{s}ov

Faculty of Mathematics and Informatics
\end{center}

\end{document}